\documentclass[12pt]{amsart}

\usepackage{amssymb}
\usepackage{amscd}
\usepackage[all]{xypic}
\swapnumbers
\def\frak{\mathfrak}
\def\Bbb{\mathbb}
\def\Cal{\mathcal}

\newtheorem*{prop*}{Proposition}

\newtheorem*{thm*}{Theorem}

\newtheorem*{lem*}{Lemma}

\newtheorem*{kor*}{Corollary}

\newcommand{\gr}{\operatorname{gr}}

\newcommand{\Ad}{\operatorname{Ad}}

\newcommand{\id}{\operatorname{id}}

\renewcommand{\ker}{\operatorname{ker}}
\newcommand{\im}{\operatorname{im}}

\newcommand{\fg}{{\frak g}}
\newcommand{\x}{\times}
\renewcommand{\o}{\circ}

\let\ccdot\cdot
\def\cdot{\hbox to 2.5pt{\hss$\ccdot$\hss}}

\newcommand{\al}{\alpha}

\newcommand{\ep}{\epsilon}

\newcommand{\ka}{\kappa}

\newcommand{\om}{\omega}
\renewcommand{\phi}{\varphi}
\newcommand{\ph}{\varphi}
\newcommand{\ps}{\psi}
\newcommand{\si}{\sigma}

\newcommand{\Ga}{\Gamma}
\newcommand{\La}{\Lambda}
\newcommand{\Om}{\Omega}
\newcommand{\Ph}{\Phi}

\def\sideremark#1{\ifvmode\leavevmode\fi\vadjust{
\vbox to0pt{\hbox to 0pt{\hskip\hsize\hskip1em
\vbox{\hsize3cm\tiny\raggedright\pretolerance10000
\noindent #1\hfill}\hss}\vbox to8pt{\vfil}\vss}}}

\begin{document}
\title[Infinitesimal Automorphisms \& Deformations]{Infinitesimal
  Automorphisms and \\ Deformations of Parabolic Geometries}
\author{Andreas \v Cap}

\address{Institut f\"ur Mathematik, Universit\"at Wien, Nordbergstr\ss
  e 15, A--1090 Wien, Austria and International Erwin Schr\"odinger
  Institute for Mathematical Physics, Boltzmanngasse 9, A--1090 Wien,
  Austria} 

\email{Andreas.Cap@esi.ac.at} 



\begin{abstract}
We show that infinitesimal automorphisms and infinitesimal
deformations of parabolic geometries can be nicely described in terms
of the twisted de--Rham sequence associated to a certain linear
connection on the adjoint tractor bundle. For regular normal
geometries, this description can be related to the underlying geometric
structure using the machinery of BGG sequences. In the locally flat
case, this leads to a deformation complex, which generalizes the
is well know complex for locally conformally flat manifolds.

Recently, a theory of subcomplexes in BGG sequences has been
developed. This applies to certain types of torsion free parabolic
geometries including, quaternionic structures, quaternionic contact
structures and CR structures. We show that for these structures one of
the subcomplexes in the adjoint BGG sequence leads (even in the curved
case) to a complex governing deformations in the subcategory of
torsion free geometries. For quaternionic structures, this deformation
complex is elliptic.
\end{abstract}

\maketitle

\section{Introduction}\label{1}
Given a smooth manifold $M$ and a type of geometric structure, it is a
natural idea to consider the moduli space, i.e.~the space of
isomorphism classes of structures of the given type on $M$. This
moduli space can be viewed as the quotient of the space of all
structures of the given type by the action of the diffeomorphism group
of $M$, which acts by pulling back structures. In general, the moduli
space is a highly complicated object. Trying to understand the moduli
space locally, one is led to the study of deformations of geometric
structures. Here deformations coming from the action of one--parameter
groups of diffeomorphisms have to be considered as trivial. Reducing
further to the formal infinitesimal level, one arrives at infinitesimal
deformations. These describe the possible directions into which a given
structure can be deformed. As before, there is the notion of a trivial
infinitesimal deformation, and the quotient of the space of all
infinitesimal deformations by the trivial ones is usually referred to
as the \textit{formal tangent space} of the moduli space at the given
structure. 

In this paper, we study infinitesimal deformations and the closely
related infinitesimal automorphisms for parabolic geometries. These
form a large class of geometric structures containing examples like
conformal, quaternionic, hypersurface type CR, and certain higher
codimension CR structures. For some of these structures, deformation
theory has been developed quite far. Infinitesimal deformations are
usually defined in an ad hoc manner as smooth sections of some
bundle. Trivial infinitesimal deformations are those which lie in the
image of some linear differential operator, whose kernel is the space
of infinitesimal automorphisms. In particular, the formal tangent
space is usually infinite dimensional. 

It is a highly interesting problem to restrict the class of allowed
deformations in such a way that one obtains a finite dimensional
moduli space. This can be done by imposing integrability conditions on
the geometric structure and looking only at deformations in the
subclass of geometries defined in that way. For parabolic geometries,
the simplest possible condition is local flatness, but in some cases
much more subtle integrability conditions can be used, for example
anti self duality for conformal structures in dimension four. 

The unifying feature of parabolic geometries is that they can be
viewed as Cartan geometries with homogeneous model a generalized flag
manifold. Regular normal geometries of this type are then equivalent
to underlying geometric structures including the examples listed above.
For Cartan geometries, there are evident notions of infinitesimal
deformations and infinitesimal automorphisms. These can be nicely
formulated in terms of a certain linear connection (which surprisingly
is different from the canonical normal tractor connection) on the
adjoint tractor bundle, see Proposition \ref{6.2}. In particular, the
relevant operators are part of the twisted de--Rham sequence
associated to this linear connection.

The machinery of Bernstein--Gelfand--Gelfand sequences (or
BGG--sequences), which was introduced in \cite{CSSBGG} and improved in
\cite{CD}, can be applied to this twisted de--Rham sequence to obtain
a sequence of higher order operators acting on sections of bundles
that can be easily interpreted in terms of the underlying structure.
For regular normal geometries, the first operator in this sequence has
the space of infinitesimal automorphisms as its kernel and the formal
tangent space to the moduli space of normal geometries as its
cokernel, see \ref{6.4} and \ref{6.5}.

For locally flat parabolic geometries (which are automatically regular
and normal, and locally isomorphic to the homogeneous model), the
twisted de--Rham sequence is a complex. Thus also the corresponding
BGG sequence is a complex which can be naturally interpreted as a
deformation complex in the category of locally flat structures.

Finally we move to more subtle integrability conditions. The recent
joint work \cite{Cap-Soucek} with V.~Sou\v cek contains a theory of
subcomplexes in BGG sequences. In that paper, we study several
examples, in which there is an interesting notion of semi flatness
which includes (and in most cases is equivalent to) torsion freeness.
In particular, these include quaternionic structures and CR
structures, but also quaternionic contact structures (torsion free
ones in dimension $7$) as introduced in \cite{Biq,Biq2}. In section
\ref{7}, we show that for all these geometries a certain subcomplex of
the adjoint BGG sequence can be naturally interpreted as a deformation
complex in the subcategory of semi flat geometries. For quaternionic
structures, this deformation complex is elliptic.

\noindent{\bf Acknowledgments.}  This work was supported by project
P15747--N05 of the ``Fonds zur F\"orderung der wissenschaftlichen
Forschung'' (FWF). Discussions with D.~Calderbank, M.~Eastwood, and
R.~Gover were very helpful. This paper grew out of the joint work
\cite{Cap-Soucek} with V.~Sou\v cek and I particularly want to thank
him for many discussions.

\section{Some background}\label{2}

We very briefly review some background. Some more details can be found
in \cite{Cap-Soucek} and much more information is available in
\cite{CSSBGG,Weyl,book}.

\subsection{Parabolic geometries}\label{2.1}
The basic data needed to define a parabolic geometry is a semisimple
Lie algebra $\fg$ endowed with a $|k|$--grading
$\fg=\fg_{-k}\oplus\dots\oplus\fg_k$ and a group $G$ with Lie algebra
$\fg$. The subgroup $P\subset G$ consisting of all elements $g\in G$
such that $\Ad(g)(\fg^i)\subset\fg^i$ for all $i$, where $\fg^i:=\fg_i
\oplus\dots\oplus\fg_k$, is a parabolic subgroup. We will also need
the subgroup $G_0\subset P$ of all elements whose adjoint action
preserves the grading of $\frak g$. 

Parabolic geometries of type $(G,P)$ are then defined as Cartan
geometries of that type. Thus such a geometry on a smooth manifold $M$
consists of a principal $P$--bundle $p:\Cal G\to M$ and a Cartan
connection $\om\in\Om^1(\Cal G,\fg)$. The \textit{homogeneous model}
of parabolic geometries of type $(G,P)$ is given by the canonical
principal bundle $G\to G/P$ with the left Maurer--Cartan form as a
Cartan connection. A morphism of parabolic geometries is a
homomorphism of principal bundles which is compatible with the Cartan
connections. In particular, any morphism is a local diffeomorphism.

The curvature of a Cartan connection $\om$ can be viewed as
$K\in\Om^2(\Cal G,\frak g)$ defined by the structure equation
$$
K(\xi,\eta)=d\om (\xi,\eta)+[\om(\xi),\om(\eta)],
$$
where $\xi$ and $\eta$ are vector fields on $\Cal G$ and the
bracket is in $\frak g$. Since $K$ is horizontal and equivariant, it
can be interpreted as a two--form $\ka$ on $M$ with values in the
associated bundle $\Cal AM:=\Cal G\x_P\frak g$, see \ref{6.1} for more
details. The bundle $\Cal AM$ is called the \textit{adjoint tractor
  bundle}.  The $P$--invariant filtration $\{\frak g^i\}$ of $\frak g$
gives rise to a filtration $\Cal AM=\Cal
A^{-k}M\supset\dots\supset\Cal A^kM$ by smooth subbundles and the Lie
bracket on $\frak g$ gives rise to a tensorial bracket $\{\ ,\ \}$ on
$\Cal AM$, making it into a bundle of filtered Lie algebras modeled on
$\frak g$.

On the other hand, the Cartan connection $\om$ induces an isomorphism
between the tangent bundle $TM$ and the associated bundle $\Cal
G\x_P(\frak g/\frak p)$. Hence there is a natural projection $\Pi:\Cal
AM\to\Cal TM$ which induces an isomorphism $\Cal AM/\Cal A^0M\cong
TM$. Via this isomorphism, the filtration of $\Cal AM$ descends to a
filtration $TM=T^{-k}M\supset\dots\supset T^{-1}M$ of the tangent
bundle by smooth subbundles. 

Applying the projection $\Pi$ to the values of $\ka$, we obtain a
$TM$--valued two--form $\ka_-$, which is called the torsion of the
Cartan connection $\om$. The geometry is called \textit{torsion free}
if this torsion vanishes.

Via the filtrations of $TM$ and $\Cal AM$, one has a natural notion of
homogeneity for $\Cal AM$--valued differential forms. In particular,
we say that $\ka$ is homogeneous of degree $\geq\ell$, if
$\ka(T^iM,T^jM)\subset\Cal A^{i+j+\ell}M$ for all $i,j=-k,\dots,-1$. A
parabolic geometry is called \textit{regular} if its curvature is
homogeneous of degree $\geq 1$. Note that torsion free parabolic
geometries are automatically regular.

For parabolic geometries, there is a uniform normalization condition.
This comes from the Kostant codifferential, which is the differential
$\partial^*:\La^k\frak p_+\otimes\fg\to\La^{k-1}\frak p_+\otimes\fg$
in the standard complex computing Lie algebra homology of $\frak
p_+:=\frak g_1\oplus\dots\oplus\frak g_k$ with coefficients in the
representation $\fg$. Now $\frak p_+$ is dual to $\frak g/\frak p$ as
a $P$--module via the Killing form, so $\Cal G\x_P(\La^k\frak
p_+\otimes\fg)\cong\La^kT^*M\otimes\Cal AM$. Since $\partial^*$ is
$P$--equivariant it induces a bundle map $\La^kT^*M\otimes\Cal
AM\to\La^{k-1}T^*M\otimes\Cal AM$ as well as a tensorial operator
$\Om^k(M,\Cal AM)\to\Om^{k-1}(M,\Cal AM)$, which we all denote by
$\partial^*$. A parabolic geometry is called \textit{normal} if and
only if $\partial^*(\ka)=0$.

Several important geometric structures like conformal structures,
almost quaternionic structures, non--degenerate CR structures of
hypersurface type, and quaternionic contact structures admit a unique
regular normal Cartan connection of type $(G,P)$ for an appropriate
choice of $(G,P)$. Usually, the underlying structure is easily encoded
into a principal $G_0$--bundle endowed with certain partially defined
differential forms. Using quite involved prolongations procedures (see
\cite{Tan, Mor, CS}) one extends this bundle to a principal
$P$--bundle and constructs a canonical regular normal Cartan
connection. This leads to an equivalence of categories between regular
normal parabolic geometries and the underlying structures. Thus
parabolic geometries offer a powerful general machinery to study a
variety of geometric structures. 

\subsection{Bernstein--Gelfand--Gelfand sequences}\label{2.2}
BGG sequences generalize the BGG resolutions of representation theory
to sequences of invariant differential operators on parabolic
geometries. They were introduced in \cite{CSSBGG} and the construction
was improved in \cite{CD}. We will briefly sketch this improved
construction for regular geometries in the special case of the adjoint
tractor bundle, more details can be also found in
\cite{Cap-Soucek,C-tw}.

The Cartan connection $\om$ induces a natural linear connection
$\nabla$, called the adjoint tractor connection, on the adjoint
tractor bundle $\Cal AM$. This in turn induced the covariant exterior
derivative 
$$
d^\nabla:\Om^k(M,\Cal AM)\to\Om^{k+1}(M,\Cal AM).
$$
The BGG machinery relates $d^\nabla$ to higher order operators
acting on sections of certain subquotient bundles. Let
$\partial^*:\La^kT^*M\otimes\Cal AM\to\La^{k-1}T^*M\otimes\Cal AM$
denote the bundle maps induced by the Kostant co\-dif\-feren\-tial. The
kernels and images of these bundle maps are natural subbundles, so we
can look at the quotient bundles $\ker(\partial^*)/\im(\partial^*)$.
By construction, they are associated to the representations $H_k(\frak
p_+,\fg)$. It turns out that the latter representations are always
completely reducible and they are algorithmically computable using
Kostant's version of the Bott--Borel--Weil theorem. Since the
associated bundles can be viewed as the fiber--wise homology groups of
the bundle $T^*M$ of Lie algebras with coefficients in the bundle
$\Cal AM$, we denote them by $H_k(T^*M,\Cal AM)$.Note that by
construction there is a natural bundle map $\pi_H:\ker(\partial^*)\to
H_k(T^*M,\Cal AM)$, and we will denote by the same symbol the
induced tensorial operator on sections.

For a normal parabolic geometry, the Cartan curvature $\ka$ by
definition is a section of $\ker(\partial^*)$, so we obtain the
section $\ka_H=\pi_H(\ka)$ of the bundle $H_2(T^*M,\Cal AM)$, which is
called the \textit{harmonic curvature}. This is a much simpler object
than $\ka$, but still a complete obstruction to local flatness. The
components of $\ka_H$ (according to the decomposition of $H_2(\frak
p_+,\frak g)$ into irreducibles) are the fundamental invariants of a
regular normal parabolic geometry. 

We have observed in \ref{2.1} that there is a natural notion of
homogeneity for $\Cal AM$--valued forms. The operators $\partial^*$
preserve homogeneities, i.e.~if $\ph\in\Om^k(M,\Cal AM)$ is
homogeneous of degree $\geq\ell$, then so is
$\partial^*(\ph)\in\Om^{k-1}(M,\Cal AM)$. For regular parabolic
geometries, also $d^\nabla$ is compatible with homogeneities. Now the
operator $\partial^*\o d^\nabla$ evidently preserves the subspace
$\Ga(\im(\partial^*))\subset\Om^k(M,\Cal AM)$. To get the machinery
going, one only needs the fact that the lowest homogeneous component
of the restriction of $\partial^*\o d^\nabla$ to
$\Ga(\im(\partial^*))$ is tensorial and invertible. Using this one
shows that the whole operator $\partial^*\o d^\nabla$ is invertible on
$\Ga(\im(\partial^*))$ and the inverse is a (by construction natural)
differential operator.

Using this inverse, one constructs a natural differential operator
$L:\Ga(H_k(T^*M,\Cal AM))\to\Om^k(M,\Cal AM)$ which is characterized
by the properties that for $\al\in\Ga(H_k(T^*M,\Cal AM))$ one has
$\partial^*(L(\al))=0$, $\pi_H(L(\al))=\al$, and $\partial^*(d^\nabla
L(\al))=0$. The first two properties say that $L$ is a differential
splitting of the tensorial projection
$\pi_H:\Ga(\ker(\partial^*))\to\Ga(H_k(T^*M,\Cal AM))$. Therefore, the
operators $L$ are referred to as the \textit{splitting operators}. The
last property implies that we can define invariant differential
operators by
$$
D:=\pi_H\o d^\nabla\o L:\Ga(H_k(T^*M,\Cal
AM))\to\Ga(H_{k+1}(T^*M,\Cal AM)),
$$
and these operators form the adjoint BGG sequence. Each of the
bundles $H_k(T^*M,\Cal AM)$ splits into a direct sum of subbundles
according to the splitting of the representation $H_*(\frak p_+,\frak
g)$ into irreducible components. Doing this in all degrees, one
obtains a pattern of operators acting between the various components.

\subsection{Infinitesimal deformations of conformal
  structures}\label{2.3} 

For the convenience of the reader, we briefly review some basic
results on infinitesimal deformations of conformal structures. Let $M$
be a smooth manifold of dimension $n\geq 3$ and let $[g]$ be a
conformal class of pseudo--Riemannian metrics on $M$. An infinitesimal
deformation of a pseudo--Riemannian metric is simply a smooth section
$h$ of the bundle $S^2T^*M$. To obtain a deformation on the conformal
class $[g]$ one first requires $h$ to be trace free, and second one
needs that rescaling the metric $g$ in the conformal class, $h$ has to
rescale in the same way. This means that $h$ has to be a section of
the tensor product of $S^2_0T^*M$ with a certain density bundle. Using
the notation and conventions of \cite{confamb}, the right bundle is
$F_1:=S^2_0T^*M[2]=S^2_0T^*M\otimes\Cal E[2]$.

Trivial deformations are those coming from pulling back the given
structure along diffeomorphisms. Viewing the conformal class $[g]$ as
a section $\mathbf{g}$ of $F_1$, this means that trivial infinitesimal
deformations are those of the form $\Cal L_\xi\mathbf{g}$, for some
vector field $\xi$ on $M$. Here $\Cal L$ denotes the Lie derivative.
Hence the quotient of all infinitesimal deformations by the trivial
ones can be interpreted as the cokernel of the (by construction
invariant) linear differential operator $D_0:\Ga(TM)\to \Ga(F_1)$
given by $D_0(\xi)=\Cal L_\xi\mathbf{g}$. It is easy to verify that
$D_0$ is the conformal Killing operator. In particular, its kernel is
the space of conformal Killing fields, i.e.~of infinitesimal conformal
isometries of $(M,[g])$.

This is about how far one can get for general conformal structures. To
proceed further one can impose some integrability condition on the
conformal structure and look at deformations in the subclass of
structures satisfying this condition. The simplest choice of such a
condition is local conformal flatness. Since this is equivalent to
vanishing of the Weyl curvature, it is natural to consider the bundle
$F_2$, in which the Weyl curvature has its values, and the operator
$D_1:\Ga(F_1)\to \Ga(F_2)$, which computes the infinitesimal change of
the Weyl curvature caused by an infinitesimal deformation of the
conformal structure. If $(M,[g])$ is locally conformally flat, then
sections in the kernel of $D_1$ correspond to infinitesimal
deformations in the subcategory of locally conformally flat
structures. Moreover, $D_1\o D_0=0$ in that case, so the quotient
$\ker(D_1)/\im(D_0)$ is exactly the formal tangent space to the moduli
space of locally conformally flat structures on $M$.

It turns out that, still in the locally conformally flat case, this
extends to a fine resolution
$$
0\to \Ga(TM)\overset{D_0}{\longrightarrow} \Ga(S^2_0T^*M[2])
\overset{D_1}{\longrightarrow} \Ga(F_2)
\overset{D_2}{\longrightarrow}\dots
\overset{D_n}{\longrightarrow}\Ga(F_n)\to 0
$$ 
of the sheaf of conformal Killing fields on $M$. Constructing this
resolution by hand is fairly involved, see the book \cite{GG}.

In the case of four dimensional conformal structures, a weaker
integrability condition is available. In this case, the bundle $F_2$
splits into the direct sum $F_2^+\oplus F_2^-$ of self dual and anti
self dual parts. Accordingly, the Weyl curvature splits as $W=W^++W^-$
and correspondingly $D_1=D_1^++D_1^-$. Given an anti self dual
conformal structure, i.e.~one such that $W^+=0$, the kernel of the
operator $D_1^+$ exactly consists of infinitesimal deformations in the
subcategory of anti self dual conformal structures. It turns out that
in this case
$$
0\to\Ga(TM)\overset{D_0}{\longrightarrow}\Ga(S^2_0T^*M[2])
\overset{D_1^+}{\longrightarrow} \Ga(F_2^+)\to 0 
$$
is a complex, which is elliptic for Riemannian signature. This is
the basis of the deformation theory for anti self dual conformal
Riemannian four manifolds, see \cite{KK} and \cite{Itoh}.

\newpage

\section{Infinitesimal automorphisms and deformations}\label{6} 

\subsection{The basic setup}\label{6.1}
Fix a parabolic geometry $(p:\Cal G\to M,\om)$ of some type $(G,P)$.
By definition, the adjoint tractor bundle $\Cal AM$ is the associated
bundle $\Cal G\x_P\frak g$ corresponding to the restriction of the
adjoint representation of $G$ to $P$. Smooth sections of this bundle
are in bijective correspondence with smooth functions $f:\Cal
G\to\frak g$ such that $f(u\cdot g)=\Ad(g^{-1})(f(u))$ for all
$u\in\Cal G$ and $g\in P$. More generally, for
$k=1,\dots,\text{dim}(M)$ the space $\Om^k(M,\Cal AM)$ can be
identified with the space $\Om^k_{\text{hor}}(\Cal G,\frak g)^P$ of
$P$--equivariant, horizontal $\frak g$--valued $k$--forms on $\Cal G$.
Here $\Ph\in\Om^k(\Cal G,\frak g)$ is horizontal, if it vanishes upon
insertion of one fundamental vector field, and $P$--equivariant if
$(r^g)^*\Ph=\Ad(g^{-1})\o\Ph$ for all $g\in P$.

Explicitly, this correspondence is given as follows. For vector fields
$\xi_1,\dots,\xi_k\in\frak X(M)$, there are $P$--invariant lifts
$\tilde\xi_1,\dots,\tilde\xi_k\in\frak X(\Cal G)$. For
$\Ph\in\Om^k(\Cal G,\frak g)$ we consider the function
$\Ph(\tilde\xi_1,\dots,\tilde\xi_k):\Cal G\to\frak g$. If $\Ph$ is
horizontal and equivariant, then this function is independent of the
choice of the lifts and $P$--equivariant, so it defines a smooth
section $\ph(\xi_1,\dots,\xi_k)$ of $\Cal AM$. One immediately
verifies that this defines an element $\ph\in\Om^k(M,\Cal AM)$. Note
that this identification is independent of the Cartan connection
$\om$.

This correspondence immediately leads to a geometric interpretation of
$\Om^1(M,\Cal AM)$: Suppose that $\tilde\om\in\Om^1(\Cal G,\frak g)$
is a second Cartan connection on $\Cal G$. Then the difference
$\tilde\om-\om\in\Om^1(\Cal G,\frak g)$ is by definition horizontal
and $P$--equivariant, and thus corresponds to an element of
$\Om^1(M,\Cal AM)$. There is an obvious notion of a deformation of the
Cartan geometry $(\Cal G\to M,\om)$ as a smooth family $\om_\tau$ of
Cartan connections on $\Cal G$ parametrized by $\tau\in
(-\ep,\ep)\subset\Bbb R$ such that $\om_0=\om$. The initial direction
of this deformation is the derivative
$\frac{d}{d\tau}|_{\tau=0}\om_\tau$ of this family at $\tau=0$. By
definition, this is the limit of $\frac{1}{\tau}(\om_\tau-\om_0)$, so
it can be interpreted as $\ph\in\Om^1(M,\Cal AM)$. On the other hand,
if $\Ph\in\Om^1(\Cal G,\frak g)$ is horizontal and $P$--equivariant,
then $\om+\Ph$ is a Cartan connection provided that it restricts to a
linear isomorphism on each tangent space. Since this is an open
condition, we can view $\Om^1(M,\Cal AM)$ as the space of all
directions of deformations of the Cartan connection $\om$, i.e.~as the
space of all infinitesimal deformations of $\om$.

{}From \ref{2.1} we know that the curvature of any Cartan connection on
$\Cal G$ is naturally interpreted as an element of $\Om^2(M,\Cal
AM)$. In particular, for a deformation $\om_{\tau}$ of $\om$, the
resulting infinitesimal change of the curvature can be viewed as an
element of $\Om^2(M,\Cal AM)$.

To discuss $\Om^0(M,\Cal AM)=\Ga(\Cal AM)$ we need a second
interpretation of $C^\infty(\Cal G,\frak g)^P$. Since $\om$
trivializes $T\Cal G$, associating to a vector field $\xi$ on $\Cal G$
the function $\om\o\xi$ defines a bijection $\frak X(\Cal G)\to
C^\infty(\Cal G,\frak g)$. Equivariancy of $\om$ immediately implies
that $(\om\o\xi)\o r^g=\Ad(g^{-1})\o(\om\o\xi)$ if and only if
$(r^g)^*\xi=\xi$, so we obtain a bijection between $\Ga(\Cal AM)$ and
the space $\frak X(\Cal G)^P$ of $P$--invariant vector fields on $\Cal
G$. Notice that $P$--invariant vector fields are automatically
projectable to vector fields on $M$, and this corresponds to the
projection $\Pi:\Cal AM\to TM$ from \ref{2.1}.

A vector field $\xi\in\frak X(\Cal G)$ satisfies $(r^g)^*\xi=\xi$ if
and only if its flow commutes with $r^g$, whenever the flow is defined.
This is true for all $g\in P$ if and only if the local flows are
principal bundle automorphisms. Thus we can view the space $\Ga(\Cal
AM)$ as the space of infinitesimal principal bundle automorphisms of
the Cartan bundle $\Cal G$.

\subsection{}\label{6.2}
Given a section of $\Cal AM$, we can look at the corresponding vector
field on $\Cal G$. The local flows of this vector field are principal
bundle automorphisms, so we can use them to pull back the Cartan
connection $\om$, which locally defines a deformation of $\om$.
Deformations obtained in this way and also the corresponding
infinitesimal deformations are called \textit{trivial}. Note that
while flows may be only locally defined the corresponding
infinitesimal deformation is always defined globally.

An automorphism of the parabolic geometry $(\Cal G,\om)$ by definition
is a principal bundle automorphism $\Ph$ of $\Cal G$ such that
$\Ph^*\om=\om$. Correspondingly, an infinitesimal automorphism is a
$P$--invariant vector field $\xi$ on $\Cal G$ such that the induced
infinitesimal deformation of the Cartan connection vanishes
identically.

In studying the infinitesimal change of curvature caused by an
infinitesimal deformation of the Cartan connection, there is an
additional subtlety. For a deformation $\om_{\tau}$ of $\om=\om_0$, we
may view the curvature $\ka_{\tau}$ of $\om_\tau$ as an element of
$\Om^2(M,\Cal AM)$, and we could simply differentiate this family of
sections. However, the identification of $\La^2T^*M\otimes\Cal AM$
with the associated bundle $\Cal G\x_P\La^2\frak p_+\otimes\frak g$,
which is used to construct operators acting on the curvature, depends
on the Cartan connection. The easiest way to take this into account is
to first convert $\ka_\tau$ into an equivariant function $\Cal G\to
\La^2\frak p_+\otimes\frak g$ using $\om_\tau$. Then one takes the
derivative of this family of functions at $\tau=0$ and converts it
back to an element of $\Om^2(M,\Cal AM)$ using $\om=\om_0$.

Finally observe that using the projection $\Pi:\Cal AM\to TM$, any
section of $\Cal AM$ has an underlying vector field on $M$. In
particular, for $s\in\Ga(\Cal AM)$ we can insert $\Pi(s)$ into a
(bundle valued) differential form on $M$, and we write $i_s$ for the
corresponding insertion operator. More generally, for
$\ph\in\Om^\ell(M,\Cal AM)$ and a vector bundle $V\to M$, we obtain an
insertion operator $i_\ph:\Om^k(M,V)\to\Om^{k+\ell-1}(M,V)$. 

\begin{prop*}
Let $(\Cal G\to M,\om)$ be a parabolic geometry with curvature
$\ka\in\Om^2(M,\Cal AM)$. Let $\nabla$ be the adjoint tractor
connection, and let $d^\nabla:\Om^k(M,\Cal AM)\to\Om^{k+1}(M,\Cal AM)$
be the corresponding covariant exterior derivative. Then we have:

\noindent
(1) For $s\in\Ga(\Cal AM)$, the infinitesimal deformation of $\om$
induced by the corresponding invariant vector field is given by
$\nabla s+i_s\ka$. In particular, $s$ is an infinitesimal automorphism
if and only if $\nabla s=-i_{s}\ka$.

\noindent
(2) For an infinitesimal deformation $\ph\in\Om^1(M,\Cal AM)$ of the
Cartan connection $\om$, the induced infinitesimal change of the
curvature is given by $d^\nabla\ph-i_\ph\ka\in\Om^2(M,\Cal AM)$.  
\end{prop*}
\begin{proof}
  (1) The derivative at zero of the pullback of $\om$ by the flow of
  $\xi$ is the Lie derivative $\Cal L_{\xi}\om\in\Om^1(\Cal G,\frak
  g)$. Evaluating this on a vector field $\eta$, we obtain
  $\xi\cdot\om(\eta)-\om([\xi,\eta])$. If $\eta$ is invariant and
  $t\in\Ga(\Cal AM)$ is the corresponding section, we can express this
  in terms of the operators on adjoint tractor fields introduced in
  \cite[section 3]{Cap-Gover}: The term $\xi\cdot\om(\eta)$
  corresponds exactly to the fundamental $D$--operator or fundamental
  derivative $D_st$, while the second term is computed in
  \cite[3.6]{Cap-Gover}. Inserting this we see that $(\Cal
  L_{\xi}\om)(\eta)$ corresponds to $D_ts+\{t,s\}+\ka(s,t)$, and by
  \cite[3.5]{Cap-Gover} the first two terms add up to
  $\nabla_{\Pi(t)}s$, which implies the result.

\noindent
(2) Let $\om_{\tau}$ be a deformation of $\om$, put
$\Ph:=\frac{d}{d\tau}|_{\tau=0}\om_{\tau}\in\Om^1_{\text{hor}}(\Cal
G,\frak g)^P$, and let $\ph\in\Om^1(M,\Cal AM)$ be the corresponding
element. Viewed as $K_{\tau}\in\Om^2(\Cal G,\frak g)$, the curvature
of $\om_\tau$ is given by 
$$
K_\tau(\xi,\eta)=d\om_{\tau}(\xi,\eta)+[\om_\tau(\xi),\om_{\tau}(\eta)]. 
$$
The derivative of this expression with respect to $\tau$ at $\tau=0$
is given by
$$
d\Ph(\xi,\eta)+[\Ph(\xi),\om(\eta)]+[\om(\xi),\Ph(\eta)]. 
$$
Choose $\xi$ and $\eta$ to be $P$--invariant and denote by $s$ and
$t$ the corresponding sections of $\Cal AM$. Inserting the definition
of the exterior derivative, we can rewrite the above as 
$$
D_s(\ph(\Pi(t)))-D_t(\ph(\Pi(s)))-\ph(\Pi([s,t]))-\{t,\ph(\Pi(s))\}+
\{s,\ph(\Pi(t))\}.
$$ 
As above, the first and last term adds up to
$\nabla_{\Pi(s)}(\ph(\Pi(t)))$ and similarly for the second and forth
term. Since $\Pi([s,t])=[\Pi(s),\Pi(t)]$, we see that
$\frac{d}{d\tau}|_{\tau=0}K_\tau$ is represented by the covariant
exterior derivative $d^\nabla\ph$.

As discussed above, we should however first convert $K_{\tau}$ into a
function using $\om_{\tau}$, which means looking at
$K_\tau(\om_\tau^{-1}(X),\om_\tau^{-1}(Y))$ for $X,Y\in\frak g$,
differentiate, and then convert the result back into a form using
$\om$. Differentiating the equation $X=\om_\tau(\om_\tau^{-1}(X))$ we
see that
$$
\tfrac{d}{d\tau}|_{\tau=0}\om_\tau^{-1}(X)=-\om^{-1}(\Ph(\om^{-1}(X))).
$$
To get the expression for the change of the curvature, we thus have to
add to $d^\nabla\ph$ the terms 
$$
-K_0\bigg(\om^{-1}(\Ph(\om^{-1}(X))),\om^{-1}(Y)\bigg)-
K_0\bigg(\om^{-1}(X),\om^{-1}(\Ph(\om^{-1}(Y)))\bigg),
$$
which exactly represent $-i_\ph\ka$. 
\end{proof}

\subsection*{Remark}
We consider infinitesimal automorphisms and deformations on the level
of the total space of the Cartan bundle here. As discussed in
\ref{2.1}, regular normal parabolic geometries are equivalent to
underlying structures. For several of these structures, notions of
infinitesimal automorphisms and deformations are available in the
literature, see \ref{2.3} for a sketch of the conformal case.  

For infinitesimal automorphisms, it is easy to see that the two
concepts are equivalent: The construction of the canonical normal
Cartan connection induces an equivalence of categories between regular
normal parabolic geometries and underlying structures. An automorphism
of the underlying structure uniquely lifts to an automorphism of the
parabolic geometry, and conversely any automorphism of a parabolic
geometry induces an automorphism of the underlying structure on the
base. Applying this to local flows of vector fields, one immediately
concludes that there is a bijective correspondence between
infinitesimal automorphisms in the two senses.  We shall see below
that this correspondence is implemented by the machinery of BGG
sequences.

In the case of infinitesimal deformations the question is a bit more
subtle, but the concepts still coincide in all cases that I am aware
of. The basic point here is the following: The underlying structures
of parabolic geometries can all be encoded as infinitesimal flag
structures, see \cite{Weyl}. These are principal $G_0$--bundles
endowed with certain partially defined differential forms. A small
deformation of the underlying structure cannot change the isomorphism
type of the principal bundle, so it can be viewed as a deformation of
the partially defined differential forms. Since the subgroup
$P_+\subset P$ is always contractible, the total space of the Cartan
bundle must be a trivial $P_+$--principal bundles over the underlying
$G_0$--bundle. Making choices, one can extend the partially defined
differential forms from above to a Cartan connection of the principal
$P$---bundle, and this transforms smooth families to smooth families.
The canonical Cartan connection can then be constructed by a
normalization process which again maps smooth families to smooth
families. This construction will be described in detail in
\cite{book}. In this way, any deformation of the underlying structure
gives rise to a deformation of the parabolic geometry, and since the
converse direction is obvious, this establishes the equivalence of the
two notions. We shall see below in examples that this correspondence
is implemented by the BGG machinery.

\subsection{A variant of the adjoint BGG sequence}\label{6.3}
Proposition \ref{6.2} suggests considering the linear connection
$\tilde\nabla$ on the bundle $\Cal AM$ which is defined by
$\tilde\nabla s=\nabla s+i_s\ka$:
\begin{lem*}
(1) For $\ph\in\Om^k(M,\Cal AM)$ we have
$d^{\tilde\nabla}\ph=d^\nabla\ph+(-1)^ki_{\ph}\ka$. 

\noindent
(2) The curvature $\tilde R$ of $\tilde\nabla$ is given by $\tilde
    R(\xi,\eta)(s)=(D_s\ka)(\xi,\eta)$, where $D_s$ denotes the
    fundamental derivative. 
\end{lem*}
\begin{proof}
(1) is a straightforward computation using the standard formula
\begin{align*}
(d^{\tilde\nabla}\ph)&(\xi_0,\dots,\xi_k)=
\sum_i(-1)^i\tilde\nabla_{\xi_i}(\ph(\xi_0,\dots,\hat
i,\dots,\xi_k))\\
&+\sum_{i<j}(-1)^{i+j}\ph([\xi_i,\xi_j],\xi_0,\dots,\hat i,\dots,\hat
j,\dots,\xi_k)
\end{align*}
for the covariant exterior derivative. 

\noindent
(2) The action of $\tilde R$ on $s\in\Ga(\Cal AM)$ can be computed as
    $d^{\tilde\nabla}\tilde\nabla s$. Inserting the definition of
    $\tilde\nabla$ and using (1), this equals 
$d^\nabla\nabla s+d^\nabla(i_s\ka)-i_{\tilde\nabla s}\ka$. The first
    term gives the action $\ka\bullet s$ of the curvature of $\nabla$,
    i.e.~$(\ka\bullet s)(\xi,\eta)=\{\ka(\xi,\eta),s\}$. Since $\ka$
    is the curvature of $\nabla$, the Bianchi identity for linear
    connections implies that $0=d^\nabla\ka$. Taking
    $t_1,t_2\in\Ga(\Cal AM)$ and expanding $0=d^\nabla\ka(t_1,s,t_2)$
    we obtain the formula 
$$
d^\nabla(i_s\ka)(t_1,t_2)=\nabla_s(\ka(t_1,t_2))
-\ka([s,t_1],t_2)-\ka(t_1,[s,t_2]).
$$
By \cite[Proposition 3.2]{Cap-Gover} we get
$\nabla_s(\ka(t_1,t_2))=D_s(\ka(t_1,t_2))+\{s,\ka(\xi,\eta)\}$. On the
other hand, \cite[Proposition 3.6]{Cap-Gover} reads as
$[s,t_1]=D_st_1-\tilde\nabla_{t_1}s$. Inserting all these facts into
the above formula for $d^{\tilde\nabla}\tilde\nabla s$, the claim
follows. 
\end{proof}

Using part (1), we conclude from Proposition \ref{6.2} that the
infinitesimal change of curvature caused by an infinitesimal
deformation of a Cartan connection is computed by
$d^{\tilde\nabla}$. 

Now suppose that we are dealing with a regular parabolic geometry
$(p:\Cal G\to M,\om)$. By definition, this means that $\ka$ is
homogeneous of degree $\geq 1$, i.e.~for $\xi\in\Ga(T^iM)$ and
$\eta\in\Ga(T^jM)$, we have $\ka(\xi,\eta)\in\Ga(\Cal
A^{i+j+1}M)$. If $\ph\in\Om^k(M,\Cal AM)$ is homogeneous of degree
$\geq\ell$, this immediately implies that $i_{\ph}\ka$ is homogeneous
of degree $\geq \ell+1$. Therefore $d^{\tilde\nabla}\ph$ is congruent
to $d^\nabla\ph$ modulo elements which are homogeneous of degree
$\geq\ell+1$. Hence the lowest possibly nonzero homogeneous components of
$d^\nabla\ph$ and of $d^{\tilde\nabla}\ph$ coincide. As pointed out in
\ref{2.2}, this is all we need to apply the BGG machinery to the
twisted de--Rham sequence induced by $\tilde\nabla$. 

We write $\tilde L:\Ga(H_k(T^*M,\Cal AM))\to\Om^k(M,\Cal AM)$ for the
splitting operators obtained by this construction. Their values
$\tilde L(\al)$ are characterized by $\partial^*(\tilde L(\al))=0$,
$\pi_H(\tilde L(\al))=\al$, and $\partial^*(d^{\tilde\nabla}\tilde
L(\al))=0$. The induced BGG operators $\tilde D_k:\Ga(H_k(T^*M,\Cal
AM))\to \Ga(H_{k+1}(T^*M,\Cal AM)))$ are given by $\tilde D_k=\pi_H\o
d^{\tilde\nabla}\o \tilde L$.

\subsection{Infinitesimal automorphisms}\label{6.4}
It is easy to relate the BGG sequence obtained from $d^{\tilde\nabla}$
to infinitesimal automorphisms:
\begin{thm*}
Let $(p:\Cal G\to M,\om)$ be a regular normal parabolic geometry of
type $(G,P)$ corresponding to a $|k|$--grading of $\frak g$. Then the
bundle $H_0(T^*M,\Cal AM)$ equals $\Cal AM/\Cal A^{-k+1}M\cong
TM/T^{-k+1}M$. The algebraic projection $\pi_H$ and the differential
operator $\tilde L$ restrict to inverse bijections between
infinitesimal automorphisms of $(p:\Cal G\to M,\om)$ and smooth
sections $\si\in\Ga(TM/T^{-k+1}M)$ such that $\tilde D_0(\si)=0$.
\end{thm*}
\begin{proof}
The bundle $H_0(T^*M,\Cal AM)$ corresponds to the representation
$H_0(\frak p_+,\frak g)$. By definition, this homology group is $\frak
g/[\frak p_+,\frak g]$, and it is well known that $[\frak p_+,\frak
g]=\fg^{-k+1}$, so the statement about $H_0(T^*M,\Cal AM)$ follows. 

By part (1) of Proposition \ref{6.2}, a smooth section $s\in\Ga(\Cal
AM)$ defines an infinitesimal automorphism if and only if
$\tilde\nabla s=0$. If this is the case, then in particular
$\partial^*(\tilde\nabla s)=0$, and since $\partial^*(s)=0$ is
automatically satisfied, this implies $s=\tilde L(\pi_H(s))$ and
$\tilde D_0(\pi_H(s))=0$. Hence $\pi_H$ restricts to an injection from
infinitesimal automorphisms to $\ker(\tilde D_0)$.

Conversely, if $\si\in\Ga(TM/T^{-k+1}M)$ satisfies $\tilde
D_0(\si)=0$, then put $s:=\tilde L(\si)$. Then
$\partial^*(\tilde\nabla s)=0$ and $\tilde D_0(\si)=0$ implies that
$\pi_H(\tilde\nabla s)=0$, so $\tilde\nabla s$ is a section of the
subbundle $\im(\partial^*)\subset T^*M\otimes\Cal AM$. By part (2) of
Proposition \ref{6.2}, we get $d^{\tilde\nabla}\tilde\nabla s=D_s\ka$
and by naturality of the fundamental derivative and normality we get
$\partial^*(D_s\ka)=D_s\partial^*(\ka)=0$. But from \ref{2.2} we know
that $\partial^*\o d^{\tilde\nabla}$ is injective on sections of
$\im(\partial^*)$, so $\tilde\nabla s=0$ and $s$ is an infinitesimal
automorphism.
\end{proof}

\subsection{}\label{6.4a}
To complete the discussion of infinitesimal automorphisms, it remains
to compare the first operator $\tilde D_0$ in the BGG sequence
associated to $\tilde\nabla$ with the first operator $D_0$ in the BGG
sequence associated to $\nabla$.

\begin{thm*}
  Let $(p:\Cal G\to M,\om)$ be a regular normal parabolic geometry of
  type $(G,P)$, and let $\frak g$ be the Lie algebra of $G$. Let $L$
  and $\tilde L$ be the splitting operators in degree zero and $D_0$
  and $\tilde D_0$ the BGG operators obtained from $\nabla$ and
  $\tilde\nabla$, respectively.

\noindent
(1) If $\frak g$ is $|1|$--graded or $(p:\Cal G\to M,\om)$ is torsion
    free, then $L=\tilde L:\Ga(TM/T^{-k+1}M)\to\Ga(\Cal AM)$ and
    $\tilde D_0(\si)=D_0(\si)+\pi_H(i_{L(\si)}\ka)$. 

\noindent
(2) If $(p:\Cal G\to M,\om)$ is torsion free and $H_1(\frak p_+,\frak
g)$ is concentrated in non--positive homogeneous degrees then $\tilde
D_0=D_0$.
\end{thm*}
\begin{proof}
  (1) We start by computing $\partial^*(i_\xi\ka)$ for an arbitrary
  vector field $\xi\in\frak X(M)$. Locally, we can write $\ka$ as a
  finite sum of terms of the form $\ph\wedge\ps\otimes t$ for
  $\ph,\ps\in\Om^1(M)$ and $t\in\Ga(\Cal AM)$. By definition,
  $\partial^*(\ka)$ is then the sum of the corresponding terms of the
  form
$$
-\ps\otimes\{\ph,t\}+\ph\otimes\{\ps,t\}-\{\ph,\ps\}\otimes t.
$$
On the other hand, $i_\xi\ka$ is the sum of the terms
$\ph(\xi)\ps\otimes t-\ps(\xi)\ph\otimes t$. Thus
$\partial^*(i_\xi\ka)$ is the sum of the terms
$\ph(\xi)\{\ps,t\}-\ps(\xi)\{\ph,t\}$, and we conclude that
$$
\partial^*(i_\xi\ka)=-i_\xi\bigg(\partial^*(\ka)-(\{\ ,\ \}\otimes\id)(\ka)\bigg),
$$
where in the last term we use $\{\ ,\ 
\}\otimes\id:\La^2T^*M\otimes\Cal AM\to T^*M\otimes\Cal AM$.  Since we
are dealing with a normal parabolic geometry, we have
$\partial^*(\ka)=0$. In the case of a $|1|$--grading the map $\{\ ,\ 
\}:\La^2T^*M\to T^*M$ is identically zero, so we get
$\partial^*(i_\xi\ka)=0$ in this case.

In the torsion free case, we fist observe that the kernel of $[\ ,\
]\otimes\id$ is a $P$--submodule in $\La^2\frak p_+\otimes\frak
g$. For any normal parabolic geometry, the harmonic curvature
$\ka_H=\pi_H(\ka)$ has values in $H_2(T^*M,\Cal AM)$. By Kostant's
version of the Bott--Borel--Weyl Theorem (see \cite{Kostant}) the
corresponding subrepresentation has multiplicity one in $\La^*\frak
p_+\otimes\frak g$. In particular, it has to be contained in the
kernel of $[\ ,\ ]\otimes\id$. By \cite[Theorem 3.2 (1)]{C-tw} the
curvature of any torsion free parabolic geometry therefore has values
in the kernel of $\{\ ,\ \}\otimes\id$, so we again conclude that
$\partial^*(i_{\xi}\ka)=0$ for each $\xi$.

For a section $\si$ of $TM/T^{-k+1}M$, consider $L(\si)$. By
construction this satisfies $\partial^*(L(\si))=0$,
$\pi_H(L(\si))=\si$, and $\partial^*(\nabla L(\si))=0$. Since
$\tilde\nabla L(\si)=\nabla L(\si)+i_{L(\si)}\ka$, so we also have
$\partial^*(\tilde\nabla L(\si))=0$. Hence $L(\si)$ satisfies the three
properties which characterize $\tilde L(\si)$ and $\tilde L=L$
follows. Using this we obtain
$$
\tilde D_0(\si)=\pi_H(\tilde\nabla
L(\si))=D_0(\si)+\pi_H(i_{L(\si)}\ka). 
$$

\noindent
(2) Since we are dealing with a torsion free geometry, we get
$i_s\ka\in\Om^1(M,\Cal A^0M)\subset\Om^1(M,\Cal AM)$ for each
$s\in\Ga(\Cal AM)$. In particular, $i_s\ka$ is always homogeneous of
degree $\geq 1$, so by the assumption on $H_1(\frak p_+,\frak g)$ we
get $\pi_H(i_{L(\si)}\ka)=0$ for any section $\si$ of $TM/T^{-k+1}M$.
\end{proof}

\begin{kor*} 
  Suppose that $(p:\Cal G\to M,\om)$ is torsion free and $H_1(\frak
  p_+,\fg)$ is concentrated in non--positive homogeneous degrees, and
  that $s\in\Ga(\Cal AM)$ satisfies $\nabla s=0$. Then $i_s\ka=0$ and
  in particular $s$ is an infinitesimal automorphism.
\end{kor*}
\begin{proof}
  Since $\nabla s=0$ we get $s=L(\pi_H(s))$ and $D_0(\pi_H(s))=0$. By
  the Theorem, we have $L=\tilde L$ and $D_0=\tilde D_0$, and in the
  proof of Theorem \ref{6.4}, we have seen that $\tilde
  D_0(\pi_H(s))=0$ implies $\tilde\nabla s=0$.
\end{proof}

\subsection*{Remark}
(1) The condition that $H_1(\frak p_+,\frak g)$ is concentrated in
non--positive homogeneous degrees is easy to verify, see \cite{Yam} or
\cite[Proposition 2.7]{CS}: The semisimple $|k|$--graded Lie algebra
$\frak g$ decomposes as a direct sum of $|k_i|$--graded simple ideals
with $k_i\leq k$ for each $i$. The condition is equivalent to the fact
that none of these simple ideals is of type $A_\ell$ or $C_\ell$ with
the grading corresponding to the first simple root. If $\frak g$
itself is simple, then this exactly excludes classical projective
structures and a contact analog of these. Note that in the latter two
cases regular normal parabolic geometries are automatically torsion
free, so part (1) holds for all regular normal geometries in these
cases.

\noindent
(2) The statement of the corollary is rather surprising even in
special cases like conformal structures. The identities responsible
for its validity are contained in the proof of Lemma \ref{6.3}. From
this proof one easily deduces $d^\nabla(i_s\ka)=D_s\ka-\ka\bullet
s+i_{\tilde\nabla s}\ka$ for any $s\in\Ga(\Cal AM)$. If $\nabla s=0$,
then $0=d^\nabla(\nabla s)=\ka\bullet s$ and if the geometry is
torsion free then this also implies that $\tilde\nabla s$ has values
in $\Cal A^0M$ and hence $i_{\tilde\nabla s}\ka=0$. Since
$0=D_s\partial^*(\ka)=\partial^*(D_s\ka)$ we obtain
$\partial^*d^\nabla(i_s\ka)=0$, which under the assumptions of the
Corollary implies $i_s\ka=0$.

\subsection{Infinitesimal deformations}\label{6.5}
We next study infinitesimal deformations of parabolic geometries.
Consider an infinitesimal deformation $\ph\in\Om^1(M,\Cal AM)$ of a
regular normal parabolic geometry $(p:\Cal G\to M,\om)$. Then $\ph$ is
called normal, if the deformed curvature (infinitesimally) remains
normal, so according to Propositions \ref{6.2} and \ref{6.3}, this is
the case if and only if $\partial^*(d^{\tilde\nabla}\ph)=0$.

The BGG machinery now easily implies that the operator $\tilde D_0$
whose kernel is the space of infinitesimal automorphisms, also has the
formal tangent space to the moduli space of normal geometries as its
cokernel:
\begin{thm*}
Let $(p:\Cal G\to M,\om)$ be a regular normal parabolic geometry.
(1) Any trivial infinitesimal deformation of $\om$ is normal.

\noindent
(2) The splitting operator $\tilde L:\Ga(H_1(T^*M,\Cal
AM))\to\Om^1(M,\Cal AM)$ induces a bijection between
$\Ga(H_1(T^*M,\Cal AM))/\im(\tilde D_0)$ and the formal tangent space
at the given structure to the moduli space of all normal parabolic
geometries on $M$.

\noindent
(3) The BGG operator $\tilde D_1$ computes the infinitesimal change of
    the harmonic curvature caused by the infinitesimal deformation
    $\tilde L(\al)$ associated to $\al\in\Ga(H_1(T^*M,\Cal AM))$.
\end{thm*}
\begin{proof}
We have already observed in the proof of Theorem \ref{6.4} that
$d^{\tilde\nabla}\tilde\nabla s=D_s\ka$ and that this has values in
the kernel of $\partial^*$, so (1) follows.

For $\al\in\Ga(H_1(T^*M,\Cal AM))$ we put $\ph:=\tilde L(\al)$. Then
by construction $\partial^*(d^{\tilde\nabla}\ph)=0$, so $\ph$ defines
a normal infinitesimal deformation. By Proposition \ref{6.2},
$d^{\tilde\nabla}\ph$ is the infinitesimal change of curvature caused
by $\ph$, and by definition $\tilde D_1(\al)=\pi_H(d^{\tilde\nabla}\ph)$,
which implies (3). 

If $\al=\tilde D_0(\si)$, then put $s=\tilde L(\si)$, so
$\al=\pi_H(\tilde\nabla s)$. Since $\partial^*(\tilde\nabla s)=0$ and
$\partial^*(d^{\tilde\nabla}\tilde\nabla s)=0$ we conclude that
$\tilde\nabla s=\tilde L(\al)$, so the resulting deformation is
trivial. Thus $\tilde L$ induces a map from the quotient
$\Ga(H_1(T^*M,\Cal AM))/\im(\tilde D_0)$ to normal infinitesimal
deformations modulo trivial infinitesimal deformations.

Suppose that $\tilde L(\al)=\tilde\nabla s$. Then in particular
$\partial^*(\tilde\nabla s)=0$, so $s=\tilde L(\pi_H(s))$. Hence
$\al=\tilde D_0(\pi_H(s))$ and our map is injective. To prove
surjectivity, suppose that $\ph\in\Om^1(M,\Cal AM)$ is any normal
infinitesimal deformation. Put $s=-\tilde Q\partial^*(\ph)$, where
$\tilde Q:\Ga(\im(\partial^*))\to\Ga(\im(\partial^*))$ is the inverse
of $\partial^*\o d^{\tilde\nabla}$, compare with \ref{2.2}. Replacing
$\ph$ by the equivalent infinitesimal deformation
$\ps=\ph+\tilde\nabla s$, we see that $\partial^*(\ps)=0$ and
$\partial^*(d^{\tilde\nabla}\ps)=0$, so $\ps=\tilde L(\pi_H(\ps))$ and
surjectivity follows.
\end{proof}

The relation between the splitting operators and the BGG operators
obtained from $d^\nabla$ respectively $d^{\tilde\nabla}$ is much more
complicated than for the first operator in the sequence. We just prove
a simple general result here which is sufficient to deal with the
cases discussed in this paper.

\begin{lem*}
  Let $(p:\Cal G\to M,\om)$ be a torsion free normal parabolic
  geometry. Suppose that $V\subset H_k(\frak p_+,\frak g)$ and
  $W\subset H_{k+1}(\frak p_+,\frak g)$ are irreducible components
  which are contained in homogeneity $\ell$ respectively $\ell+1$.
  Then the components of the BGG operators $\tilde D_k$ and $D_k$,
  which map sections of $\Cal G\x_PV$ to sections of $\Cal G\x_PW$,
  coincide.
\end{lem*}
\begin{proof}
  Consider a section $\al\in\Ga(\Cal G\x_PV)$ and put
  $\ph:=L(\al)\in\Om^k(M,\Cal AM)$. Then $\ph$ is homogeneous of
  degree $\geq\ell$, $\partial^*(\ph)=0$ and $\pi_H(\ph)=\al$. By part
  (1) of Lemma \ref{6.3} we get
  $d^{\tilde\nabla}\ph=d^\nabla\ph+(-1)^ki_\ph\ka$ and therefore
  $\partial^*(d^{\tilde\nabla}\ph)=(-1)^k\partial^*(i_\ph\ka)$. By
  torsion freeness $\ka$ is homogeneous of degree $\geq 2$, so
  $i_\ph\ka$ is homogeneous of degree $\geq\ell+2$. Denoting by
  $\tilde Q$ the operator used in the proof of the Theorem, we
  conclude that $\ps:=(-1)^{k+1}\tilde Q\partial^*(i_\ph\ka)$ is
  homogeneous of degree $\geq\ell+2$. By construction
  $\partial^*(\ph+\ps)=0$, $\pi_H(\ph+\ps)=\al$, and
  $\partial^*(d^{\tilde\nabla}(\ph+\ps))=0$, which implies $\tilde
  L(\al)=\ph+\ps$. Now
$$
d^{\tilde\nabla}(\ph+\ps)=d^\nabla\ph+(-1)^ki_\ph\ka+d^{\tilde\nabla}\ps,
$$
and the last two terms are homogeneous of degree $\geq\ell+2$. By
homogeneity, these terms cannot contribute to the component of the
image under $\pi_H$ that we are interested in.
\end{proof}

\subsection{On regularity}\label{6.5a}
To get a complete correspondence to underlying structures, one has to
single out regular normal infinitesimal deformations among all normal
ones. Here a normal infinitesimal deformation $\ph\in\Om^1(M,\Cal AM)$
is called regular if and only if $d^{\tilde\nabla}\ph\in\Om^2(M,\Cal
AM)$ is homogeneous of degree $\geq 1$. Notice that this condition is
vacuous if the geometry corresponds to a $|1|$--grading, and Theorem
\ref{6.5} therefore gives a complete description of the formal tangent
space to the moduli space of regular normal geometries.

In general, we can first show that trivial infinitesimal deformations
of regular normal geometries are regular. Indeed, from part (2) of
Lemma \ref{6.3} we know that for $s\in\Ga(\Cal AM)$ we have
$d^{\tilde\nabla}\tilde\nabla s=D_s\ka$. If we start from a regular
normal geometry, then $\ka$ is homogeneous of degree $\geq 1$, and by
naturality of the fundamental derivative the same is true for
$D_s\ka$. Theorem \ref{6.5} now directly implies
\begin{kor*}
  Let $(p:\Cal G\to M,\om)$ be a regular normal parabolic geometry.
  Then the formal tangent space at the given structure to the moduli
  space of regular normal geometries is the quotient of the space of
  all $\al\in\Ga(H_1(T^*M,\Cal AM))$ such that $d^{\tilde\nabla}\tilde
  L(\al)\in\Om^2(M,\Cal AM)$ is homogeneous of degree $\geq 1$ by the
  image of $\tilde D_0$.
\end{kor*}

For any concrete choice of structure, the condition on the homogeneity
of $d^{\tilde\nabla}\tilde L(\al)$ can be made more explicit by
projecting out step by step the lowest possibly nonzero homogeneous
components of $d^{\tilde\nabla}\o \tilde L$. For structures
correspondig to $|2|$--gradings, we can give a nicer description,
which will be useful in the examples in section \ref{7}.

\begin{prop*}
  Suppose that $P\subset G$ corresponds to a $|2|$--grading of $\frak
  g$. Then for any regular normal parabolic geometry $(p:\Cal G\to
  M,\om)$ of type $(G,P)$ and any section $\al\in\Ga(H_1(T^*M,\Cal
  AM))$ the form $d^{\tilde\nabla}\tilde L(\al)\in\Om^2(M,\Cal AM)$ is
  homogeneous of degree $\geq 1$ if and only $\tilde D_1(\al)$ is
  homogeneous of degree $\geq 1$.
\end{prop*}
\begin{proof}
By definition, we have $\tilde D_1(\al)=\pi_H(d^{\tilde\nabla}\tilde
L(\al))$. If $\tilde D_1(\al)$ is homogeneous of degree $\geq 1$, then
so is $\tilde L(\tilde D_1(\al))$, which differs from $d^{\tilde\nabla}\tilde
L(\al)$ by a section of $\im(\partial^*)$. Since we deal with a
$|2|$--grading, any element of $\La^3T^*M\otimes\Cal A^M$ is
homogeneous of degree $\geq 1$, and the result follows since
$\partial^*$ preserves homogeneities. 
\end{proof}

Since any irreducible component of $H_2(\frak p_+,\frak g)$ is
contained in some homogeneous degree, the condition in the proposition
simply means that all components of $\tilde D_1(\al)$ in bundles
corresponding to irreducible pieces in homogeneity zero have to
vanish.

\subsection{The locally flat case}\label{6.6}
As a simple consequence of Theorem \ref{6.5}, we can deal with the
case of locally flat geometries. The following result was first proved
in \cite{CD}.
\begin{thm*}
Let $(p:\Cal G\to M,\om)$ be a locally flat parabolic geometry. Then
the BGG sequence associated to the adjoint representation is a
complex. It can be naturally viewed as a deformation complex, i.e.~its
homologies in degrees zero and one are the space of infinitesimal
automorphisms respectively the formal tangent space to the moduli
space of all locally flat parabolic geometries on $M$.  
\end{thm*}
\begin{proof}
  By local flatness, $\nabla=\tilde\nabla$ and this connection is
  flat, so the twisted de-Rham sequence is a complex. This easily
  implies that $L\o D=d^\nabla\o L$, so the BGG sequence also is a
  complex. By Theorem \ref{6.4a}, the cohomology of this complex in
  degree zero is isomorphic to the space of infinitesimal
  automorphisms. For $\al\in\Ga(H_1(T^*M,\Cal AM))$ with $D_1(\al)=0$
  we have $d^\nabla L(\al)=LD(\al)=0$, so the infinitesimal
  deformation $L(\al)$ does not change the curvature infinitesimally.
  Since conversely $d^\nabla L(\al)=0$ clearly implies $D_1(\al)=0$,
  we see that the kernel of $D_1$ exactly corresponds to the
  infinitesimal deformations in the subcategory of locally flat
  geometries. Now the interpretation of the first cohomology follows
  from Theorem \ref{6.5}.
\end{proof}

\section{Deformation complexes for torsion free geometries}\label{7}
In the recent joint work \cite{Cap-Soucek} with V.~Sou\v cek, we have
developed a theory of subcomplexes in curved BGG sequences. This
theory applies to torsion free geometries of certain types. To have
interesting examples, one needs assumptions on the structure of the
homology groups $H_2(\frak p_+,\frak g)$, which form the degree two
part of the adjoint BGG sequence and governs the structure of the
harmonic curvature. The main examples of this situation are the ones
discussed in \cite{Cap-Soucek}. 

Here we find for all these examples a certain subcomplex in the
adjoint BGG sequence (obtained from $d^\nabla$). Using the results of
section \ref{6}, we can show that the first two operators in this
subcomplex coincide with their counterparts in the BGG sequence
obtained from $d^{\tilde\nabla}$. This leads to an interpretation of
the subcomplex as a deformation complex in the appropriate subcategory
of torsion free geometries.

\subsection{Grassmannian structures}\label{6.7}
An almost Grassmannian structure on a manifold $M$ of dimension $2n$
is essentially given by two auxiliary bundles $E$ and $F$ over $M$ of
rank $2$ respectively $n$, and an isomorphism $\Ph:E^*\otimes F\to
TM$. The bundles $E$ and $F$ are the basic building blocks for bundles
over $M$ corresponding to irreducible representations of $P$.

The BGG sequences in this case have triangular shape, see
\cite[3.4]{Cap-Soucek}. The bundle in degree $k$ of the BGG sequence
splits as a direct sum of irreducible subbundles $\Cal H_{p,q}$ with
$p+q=k$ and $0\leq p\leq q\leq n$. In particular, the second bundle
splits as $\Cal H_{0,2}\oplus\Cal H_{1,1}$, and correspondingly there
are two irreducible components in the harmonic curvature. Let us now
restrict to the case $n>2$, the case $n=2$ will be discussed below.
The harmonic curvature component in $\Ga(\Cal H_{0,2})$ is called
the torsion of the almost Grassmannian structure. Vanishing of this
torsion is equivalent to torsion freeness in the sense of
$G$--structures, and the corresponding geometries are called
Grassmannian rather than almost Grassmannian. The harmonic curvature
component in $\Cal H_{1,1}$ is a true curvature. It is shown in
\cite[Theorem 3.5]{Cap-Soucek} that in the case of Grassmannian
structures for any $p=0,\dots,n$ the parts $\Cal
H_{p,p}\to\dots\to\Cal H_{p,n}$ and for any $q=0,\dots,n$ the parts
$\Cal H_{0,q}\to\dots\to\Cal H_{q,q}$ are subcomplexes in each BGG
sequence.

The representations inducing the bundles in the adjoint BGG sequence
are determined in \cite[4.1]{Cap-Soucek}, where we have to take
$k=\ell=1$. For $j<n$, one obtains $\Cal H_{0,j}=(S^jE\otimes
E^*)_0\otimes (\La^jF^*\otimes F)_0$, where the subscript $0$ denotes
the tracefree part.  In particular $\Cal H_{0,0}=E^*\otimes F=TM$,
which also follows from Theorem \ref{6.4}, and $\Cal H_{0,1}=\frak
s\frak l(E)\otimes\frak s\frak l(F)$.  Evidently, $\Cal
H_{0,j}M\subset\La^j(E\otimes F^*)\otimes(E^*\otimes
F)=\La^jT^*M\otimes TM$. Looking at homogeneities, this implies that
the BGG operators $\Ga(\Cal H_{0,j-1})\to\Ga(\Cal H_{0,j})$ are first
order for all $j=1,\dots,n-1$. Finally, $\Cal H_{0,n}=(S^{n+1}E\otimes
E^*)_0\otimes\La^nF^*$, and the last BGG operator $\Ga(\Cal
H_{0,n-1})\to\Ga(\Cal H_{0,n})$ is of second order.

Finally, we need the bundle $\Cal H_{1,1}$ which turns out to be the
highest weight part in $\La^2E\otimes S^2F^*\otimes\frak s\frak
l(F)$. This is contained in $\La^2T^*M\otimes L(TM,TM)$, so the BGG
operator $\Ga(\Cal H_{0,1})\to\Ga(\Cal H_{1,1})$ is a second order
operator.

\begin{thm*}
Let $M$ be a Grassmannian manifold of dimension $2n\geq 6$. Then the
subcomplex
$$
0\to \Ga(\Cal H_{0,0})\to\Ga(\Cal H_{0,1})\to\dots\to\Ga(\Cal
H_{0,n})\to 0
$$
of the adjoint BGG sequence is a deformation complex in the
subcategory of Grassmannian structures.
\end{thm*}
\begin{proof}
  The first two operators in this sequence are just the first two
  operators in the full adjoint BGG sequence, and from Theorem
  \ref{6.4a} and Lemma \ref{6.5} we conclude that they coincide with
  their counterparts constructed from $\tilde\nabla$ rather than
  $\nabla$. The statement on the cohomology in degree zero then
  follows from Theorem \ref{6.4}.
  
  By part (2) of Theorem \ref{6.5} and since regularity is automatic
  for $|1|$--gradings, the quotient $\Ga(\Cal H_{0,1})/\im(D_0)$ is isomorphic to infinitesimal deformations
  of $M$ in the category of almost Grassmannian structures modulo
  trivial infinitesimal deformations. On the other hand, part (3) of
  Theorem \ref{6.5} implies that the kernel of $\Ga(\Cal
  H_{0,1}M)\to\Ga(\Cal H_{0,2}M)$ corresponds exactly to those
  deformations for which the infinitesimal change of torsion is
  trivial, so these are exactly the infinitesimal deformations in the
  category of Grassmannian structures.
\end{proof}

\subsection*{Remark}
(1) For almost Grassmannian structures, the right definition of an
infinitesimal deformation is not immediately evident. Is is a nice
feature of the approach via parabolic geometries and the BGG
machinery, that it shows that infinitesimal deformations are smooth
sections of the bundle $\frak s\frak l(E)\otimes\frak s\frak l(F)$.
This can be seen directly as follows.

The only part of an almost Grassmannian structure that can be deformed
nontrivially is the isomorphism $\Ph:E^*\otimes F\to TM$.
Infinitesimally, deformations of this isomorphisms are linear maps
$E^*\otimes F\to TM$ modulo those, which are compatible with $\Ph$.
Using $\Ph$ to convert the target of such a map back to $E^*\otimes
F$, these are exactly endomorphisms of $E^*\otimes F$ modulo those
which are of the form $\ph\otimes\id_F+\id_E\otimes\ps$.

\noindent
(2) By Theorem \ref{6.5}, the splitting operator $\tilde L:\Ga(\Cal
H_{0,1})\to\Om^1(M,\Cal AM)$ computes the infinitesimal deformation of
the canonical Cartan connection caused by an infinitesimal deformation
of the underlying structure.

\subsection{The case $n=2$}\label{6.8}
In this case, $\dim(M)=4$ and an almost Grassmannian structure is
equivalent to a conformal spin structure with split signature $(2,2)$.
Basically, this is due to the fact that $SL(4,\Bbb R)$ naturally is a
two fold covering of $SO(3,3)$. Here the situation is more symmetric
than for general Grassmannian structures and the two components of the
harmonic curvature are the self dual and the anti self dual parts of
the Weyl curvature. Theorem \ref{6.5} directly leads to a complex
$$
0\to\Ga(\Cal H_{0,0})\to\Ga(\Cal H_{0,1})\to\Ga(\Cal H_{0,2})\to 0
$$
inside the BGG sequence obtained from $d^{\tilde\nabla}$, and, for
anti self dual structures, an interpretation as a deformation complex
in the category of anti self dual conformal structures. This is
exactly the split signature version of the complex discussed in
\ref{2.3}. However, in this case the second operator (which has order
two) differs (tensorially) from its counterpart in the standard
adjoint BGG sequence.

\subsection{Quaternionic structures}\label{6.10}
An almost quaternionic structure on a smooth manifold of dimension
$4n$ is given by a rank 3 subbundle $Q\subset L(TM,TM)$ which is
locally spanned by three almost complex structures $I$, $J$, and
$K=IJ=-JI$. However, these local almost complex structures are an
additional choice and not an ingredient of the structure.
Equivalently, one can view an almost quaternionic structure are a
reduction of the structure group of the linear frame bundle to the
subgroup $S(GL(1,\Bbb H)GL(n,\Bbb H))\subset GL(4n,\Bbb R)$. Replacing
$S(GL(1,\Bbb H)GL(n,\Bbb H))$ by the two--fold covering $S(GL(1,\Bbb
H)\x GL(n,\Bbb H))$ one has an equivalent description as an
identification of the complexified tangent bundle $TM\otimes\Bbb C$
into the tensor product $E\otimes F$, where $E$ has complex rank two
and $F$ has complex rank $2n$. Hence after complexification we are in
the same situation as for almost Grassmannian structures with even
dimensional $F$.

In particular, the BGG sequences have the same shape as in the almost
Grassmannian case, and the operators have the same orders. Moreover,
after complexification the bundles showing up in each BGG sequence are
the same as in the almost Grassmannian case. In particular, there are
again two harmonic curvature components and for $n>1$ (the case $n=1$
will be discussed below) one of them is a torsion and the other is a
true curvature. Vanishing of the torsion is again equivalent to
torsion freeness in the sense of $G$--structures and the corresponding
geometries are referred to as quaternionic rather than almost
quaternionic. For quaternionic structures one obtains subcomplexes in
all BGG sequences which have the same form as in the Grassmannian
case. 

\begin{thm*}
Let $M$ be a quaternionic manifold of dimension $4n\geq 8$. Then the
subcomplex 
$$
0\to \Ga(\Cal H_{0,0})\to\Ga(\Cal H_{0,1})\to\dots\to\Ga(\Cal
H_{0,n})\to 0
$$
of the adjoint BGG sequence is an elliptic complex, which can be
naturally interpreted as a deformation complex in the category of
quaternionic structures.
\end{thm*}
\begin{proof}
  The interpretation as a deformation complex works exactly as in the
  Grassmannian case. Exactness of the symbol sequence is the special
  case $k=\ell=1$ of \cite[Theorem 4.3]{Cap-Soucek}.
\end{proof}

Describing the bundles which show up in the deformation complex is
straightforward, but we do not need this description here. Let us just
note the description of the bundle $\Cal H_{0,1}$, whose sections are
the infinitesimal deformations of an almost quaternionic structure.
For $Q\subset L(TM,TM)$ let $L_Q(TM,TM)\subset L(TM,TM)$ denote the
subbundle of those endomorphisms which commute with any element of
$Q$. Then it turns out that $L(TM,TM)\cong Q\oplus L_Q(TM,TM)\oplus
\Cal H_{0,1}$ and that $\Cal H_{0,1}$ is isomorphic to the tensor
product of $Q$ with the space of tracefree elements of $L_Q(TM,TM)$.

In the special case $n=1$, an almost quaternionic structure on a four
manifold is equivalent to a conformal Riemannian spin structure. As in
\ref{6.8}, we obtain the deformation complex discussed in \ref{2.3}
directly from Theorem \ref{6.5}. Again, the second operator in the
sequence differs from the one in the standard adjoint BGG sequence.
Ellipticity can be easily verified directly. 

\subsection{Lagrangean contact structures}\label{6.12}
A Lagrangean contact structure on a smooth manifold $M$ of dimension
$2n+1$ is given by a codimension one subbundle $H\subset TM$, which
defines a contact structure on $M$, together with a decomposition
$H=E\oplus F$ as a direct sum of two Lagrangean (or Legendrean)
subbundles. This means that the Lie bracket of two sections of $E$
(respectively $F$) is a section of $H$. Since $H$ defines a contact
structure, this forces $E$ and $F$ to be of rank at most (and hence
equal to) $n$. We will assume $n\geq 2$ throughout.

The form of the BGG sequences is described in \cite[3.6]{Cap-Soucek}.
For $k\leq n$, the bundle in degree $k$ of each BGG sequence splits as
$\oplus_{p,q}\Cal H_{p,q}$ with $p+q=k$ and $0\leq p,q$. The
decomposition of the bundle in degree $n+k$ has the same form as for
degree $n-k+1$. In particular, there are three components in the
harmonic curvature. The components in $\Ga(\Cal H_{2,0})$ and
$\Ga(\Cal H_{0,2})$ are torsions, which are exactly the obstructions
to integrability of the subbundles $E$ and $F$ of $TM$. The component
in $\Ga(\Cal H_{1,1})$ is a true curvature. For torsion free
geometries, the bundles $E$ and $F$ are integrable, so $M$ locally
admits two transversal fibrations onto manifolds of dimension $n+1$
such that the two vertical bundles span a contact distribution on $M$
and both are Legendrean.  In the torsion free case, there many
subcomplexes in each BGG sequence (see \cite[3.7]{Cap-Soucek}), in
particular the bundles $\Cal H_{0,q}$ for $q=0,\dots, n$ and $\Cal
H_{p,0}$ for $p=0,\dots,n$ both form subcomplexes.

Specializing to the adjoint BGG sequence, we know from \ref{6.4} that
the first bundle $\Cal H_{0,0}$ is the quotient $Q:=TM/H$ of the
tangent bundle by the contact subbundle. For $0<j\leq n$ one easily
verifies that the bundle $\Cal H_{0,j}$ can be described as follows.
Since $E$ and $F$ are Legendrean, the contact structure defines
isomorphisms $F\cong E^*\otimes Q$. Thus $\La^jE^*\otimes F\cong
\La^jE^*\otimes E^*\otimes Q$, and $\Cal H_{0,j}\subset
\La^jE^*\otimes F$ corresponds to the kernel of the alternation
$\La^jE^*\otimes E^*\otimes Q\to\La^{j+1}E^*\otimes Q$. In particular,
the operator mapping sections of $\Cal H_{0,0}$ to sections of $\Cal
H_{0,1}$ must be of second order, while for $1\leq j<n$ the operator
$\Ga(\Cal H_{0,j})\to\Ga(\Cal H_{0,j+1})$ is first order. In the same
way, the bundles $\Cal H_{j,0}$ for $1\leq j\leq n$ can be described
as subbundles in $\La^jF^*\otimes E$ and one gets the analogous
results for the orders of the operators.

\begin{thm*}
  Let $(M,E,F)$ be a torsion free Lagrangean contact structure. Then
  the subcomplex
$$
0\to\Ga(\Cal H_{0,0})\to \begin{array}c\Ga(\Cal H_{0,1})\\[-2pt] \oplus\\
\Ga(\Cal H_{1,0})\end{array}\to\dots \to
\begin{array}c\Ga(\Cal H_{0,n})\\[-2pt] \oplus\\ \Ga(\Cal H_{n,0})
\end{array}\to 0
$$
in the adjoint BGG sequence can be naturally viewed as a
deformation complex in the category of torsion free Lagrangean contact
structures.
\end{thm*}
\begin{proof}
  From Theorem \ref{6.4a} and Lemma \ref{6.5} we see that the first
  two operators in this complex coincide with their counterparts in
  the BGG sequence obtained from $d^{\tilde\nabla}$. Since the bundles
  $\Cal H_{2,0}$, $\Cal H_{1,1}$, and $\Cal H_{0,2}$ are all contained
  in positive homogeneous degrees, normal infinitesimal deformations
  are automatically regular by Proposition \ref{6.5a}. Now the
  interpretation as a deformation complex works as for Grassmannian
  structures.
\end{proof}

We can again see directly that sections of the bundle $\Cal
H_{0,1}\oplus \Cal H_{1,0}$ are the right notion for infinitesimal
deformations of a Lagrangean contact structure. Since contact
structures are rigid, the only way to deform such a structure is
deforming the decomposition $H=E\oplus F$. Infinitesimally, a
deformation of the subbundle $E\subset H$ is given by a linear map
$E\to H/E\cong F$. Such a deformation is in the direction of a
Legendrean subbundle if and only if the corresponding map $E\x E\to Q$
has trivial alternation. Likewise, $\Ga(\Cal H_{1,0})$ describes
infinitesimal Lagrangean deformations of $F\subset H$. The splitting
operator $\tilde L_1$ again computes the infinitesimal deformation of
the normal Cartan connection caused by an infinitesimal deformation of
a Lagrangean contact structure.

Note that the deformation complex cannot be elliptic or subelliptic,
since $Q$ is a real line bundle, while all other bundles showing up in
the subcomplex have even rank.

\subsection{CR structures}\label{6.13}
This case is closely parallel to the case of Lagrangean contact
structures. The geometries in question are non--de\-ge\-ne\-rate almost CR
structures of hypersurface type, which satisfy partial integrability,
a weakening of the usual integrability condition for CR structures,
see \cite[4.15]{CS}.  Compared the the case of Lagrangean contact
structures, one only replaces the decomposition of the contact
subbundle (which can be interpreted as an almost product structure) by
an almost complex structure $J$ on the contact subbundle $H$.  The
condition that the two subbundles are Legendrean corresponds to the
partial integrability condition for almost CR structures. The
complexification $H\otimes\Bbb C$ splits as $H^{1,0}\oplus H^{0,1}$
into holomorphic and anti holomorphic part, so an this level the
picture is parallel to the Lagrangean contact case.

In particular, complex BGG sequences have exactly the same form as for
Lagrangean contact structures. However, BGG sequences corresponding to
real representations without an invariant complex structure (like the
adjoint representation) are different. They are obtained by
``folding'' a complex BGG pattern, see \cite[3.8]{Cap-Soucek}. In
particular, there are only two irreducible components in the harmonic
curvature. One of these components is a torsion (corresponding to the
two torsions for Lagrangean contact structures), while the other is a
curvature. The torsion is a multiple of the Nijenhuis tensor, so the
torsion free geometries are exactly  CR structures, see
\cite[4.16]{CS}. For CR structures, one obtains many subcomplexes in
BGG sequences, see \cite[Theorem 3.8]{Cap-Soucek}.

Using the notation of \cite[3.8]{Cap-Soucek}, there is a subcomplex in
the adjoint BGG sequence which starts at $\Cal H_{0,0}$. This has the
form
$$
0\to \Ga(\Cal H_{0,0})\to\Ga(\Cal H_{1,0})\to\dots\to\Ga(\Cal
H_{n,0})\to 0. 
$$
and apart from $\Cal H_{0,0}=Q:=TM/H$, all the bundles $\Cal
H_{j,0}$ in the sequence are complex vector bundles. To identify
them, we just have to observe that their complexification splits into
a direct sum of two complex vector bundles which exactly correspond to
the two bundles in the Lagrangean contact case with $E$ and $F$
replaced by $H^{1,0}$ and $H^{0,1}$. In particular we see that $\Cal
H_{j,0}\otimes\Bbb C\subset L_{\Bbb C}(\La^j(H\otimes\Bbb
C)^*,H\otimes\Bbb C)$ and the components are singled out by their
complex (anti--)linearity properties. For example, $\Cal
H_{0,1}\otimes\Bbb C$ is contained in $L(H^{1,0},H^{0,1})\oplus
L(H^{0,1},H^{1,0})$, which exactly means that $\Cal H_{0,1}$ consists
of conjugate linear maps $H\to H$. A conjugate linear map $\ph$ lies
in $\Cal H_{0,1}$ if and only if the corresponding bilinear map $H\x
H\to Q$ is symmetric. In particular, conjugate linear maps are exactly
infinitesimal deformations of the almost complex structure $J$ (which
is the only deformable ingredient in the structure) and the symmetry
condition takes care about partial integrability.  As before we
deduce:

\begin{thm*}
  Let $(M,H,J)$ be a non--degenerate CR structure of hypersurface
  type. Then the subcomplex
  $$
  0\to \Ga(\Cal H_{0,0})\to\Ga(\Cal H_{1,0})\to\dots\to\Ga(\Cal
  H_{n,0})\to 0
$$
in the adjoint BGG sequence can be naturally interpreted as a
deformation complex in the category of CR structures.
\end{thm*}

\subsection*{Remark} 
This deformation complex has been found (by ad hoc methods) and
successfully applied to the deformation theory of strictly
pseudoconvex compact CR manifolds earlier. The part starting from
$\Cal H_{1,0}M$ is used in the work of T.~Akahori in the case $n\geq
3$, see e.g.~\cite{Akahori}. The full complex was constructed in
\cite{AGL} for $n=2$. Since the first bundle in the complex is a real
line bundle while all others a complex vector bundles, there is again
no hope for the whole complex to be elliptic or subelliptic.
Nonetheless, for some of the operators in the complex one can prove
subelliptic estimates (in the strictly pseudoconvex case), which play
a crucial role in the applications to deformation theory.

\subsection{Quaternionic contact structures}\label{6.14}
These geometries are given by certain codimension three subbundles in
the tangent bundles of manifolds of dimension $4n+3$. Recall first
that for $p+q=n$, there is (up to isomorphism) a unique quaternionic
Hermitian form of signature $(p,q)$ on $\Bbb H^n$. The imaginary part
of this form is a skew symmetric bilinear map $\Bbb H^n\x\Bbb
H^n\to\im(\Bbb H)$. Putting $\frak g_1:=\Bbb H^n$ and $\frak
g_2:=\im(\Bbb H)$ this imaginary part makes $\frak g_1\oplus\frak g_2$
into a nilpotent graded Lie algebra, called the quaternionic
Heisenberg algebra of signature $(p,q)$. Since the forms of signature
$(p,q)$ and $(q,p)$ differ only by sign, we may assume $p\geq q$.
Similarly, one may look at the algebra of split quaternions, for which
there is a unique Hermitian form in each dimension. Correspondingly,
one obtains a unique split quaternionic Heisenberg algebra of
dimension $4n+3$ for each $n\geq 1$.

Recall that if $H\subset TM$ is a subbundle in the tangent bundle of a
smooth manifold $M$, then the Lie bracket of vector fields induces a
tensorial map $\Cal L:H\x H\to TM/H=:Q$. For any $x\in M$ we put
$\gr_1(T_xM)=H_x$ and $\gr_2(T_xM):=Q_x$. Then we can view $\Cal L$
as defining on each of the spaces
$\gr(T_xM)=\gr_1(T_xM)\oplus\gr_2(T_xM)$ the structure of a nilpotent
graded Lie algebra. A \textit{quaternionic contact structure} of
signature $(p,q)$ on a smooth manifold $M$ of dimension $4(p+q)+3$ is
a smooth subbundle $H\subset TM$ of corank $3$ such that for each
$x\in M$ the nilpotent graded Lie algebra $\gr(T_xM)$ is isomorphic to
the quaternionic Heisenberg algebra of signature $(p,q)$. A
\textit{split quaternionic contact structure} on a smooth manifold of
dimension $4n+3$ is defined similarly using the split quaternionic
Heisenberg algebra.

For $n=1$ we have $\dim(M)=7$ and there is only one possible
signature. It turns out that both the quaternionic and the split
quaternionic Heisenberg algebra are rigid in this case. Moreover,
corank three distributions defining quaternionic and split
quaternionic contact structures the two generic types of rank 4
distributions in dimension $7$. In particular, a generic real
hypersurface in a two--dimensional (split) quaternionic vector space
carries a (split) quaternionic contact structure. 

For $n>1$, there are no generic distributions of rank $4n$ in
manifolds of dimension $4n+3$, but it is known from the works of
O.~Biquard, see \cite{Biq,Biq2}, that there are many examples of
quaternionic contact structures of signature $(n,0)$.

For all these structures, the BGG sequences have the same form, see
\cite[3.9]{Cap-Soucek}. For $k=0,\dots,2n+1$ the bundle in degree $k$
splits into a direct sum of bundles $\Cal H_{p,q}$ with $p+q=k$ and
$p\geq q$, and for the bundle in degree $2n+1+k$ decomposes in the
same way as the one in degree $2n+2-k$. In particular, in degree two
we obtain two irreducible components $\Cal H_{2,0}$ and $\Cal
H_{1,1}$. The harmonic curvature component having values in the bundle
$\Cal H_{2,0}$ of the adjoint BGG sequence is a torsion, while the one
having values in $\Cal H_{1,1}$ is a curvature. For $n=1$, one obtains
a subcategory of torsion free (split) quaternionic contact structures.
However, for $n>1$, bundle $\Cal H_{2,0}$ is contained in homogeneity
zero, so vanishing of the corresponding harmonic curvature component
is forced by regularity, and any (split) quaternionic contact
structure is automatically torsion free.

By \cite[Theorem 3.10]{Cap-Soucek} there is a number of subcomplexes
in each BGG sequence for a manifold endowed with a torsion free
(split) quaternionic contact structure. In particular, the bundles
$\Cal H_{p,0}$ with $p=0,\dots,2n+1$ form a subcomplex.  For the
adjoint BGG sequence, one verifies directly that the operator $\Ga(\Cal
H_{n,0})\to\Ga(\Cal H_{n+1,0})$ is of second order, while all other
operators in the subcomplex are of first order.

\begin{thm*}
  Let $M$ be a smooth manifold of dimension $4n+3\geq 11$ endowed with
  a quaternionic contact structure or split quaternionic contact
  structure. Then the subcomplex
$$
0\to \Ga(\Cal H_{0,0})\to\Ga(\Cal H_{1,0})\to\dots\to\Ga(\Cal
H_{2n+1,0})\to 0 
$$
of the adjoint BGG sequence can be naturally interpreted as a
deformation complex in the category of torsion free (split)
quaternionic contact structures. 
\end{thm*}
\begin{proof}
  By Theorem \ref{6.4a} and Lemma \ref{6.5} the first two operators in
  this sequence coincide with their counterparts obtained from
  $d^{\tilde\nabla}$. Using Proposition \ref{6.5a} we conclude that
  the kernel of the operator $\Cal H_{1,0}\to\Cal H_{2,0}$ exactly
  corresponds to regular normal deformations.  The interpretation as a
  deformation complex then works as before.
\end{proof}

\subsection*{Remark}
The situation in the seven--dimensional case is not completely clear.
The problem here is that the operator $\Cal H_{1,0}\to\Cal H_{2,0}$ in
the adjoint BGG sequence is of second order. It seems that the two
operators obtained from $d^\nabla$ respectively from
$d^{\tilde\nabla}$ differ (tensorially) from each other.  Therefore,
there seems to be no direct way to relate the BGG sequence based on
$d^\nabla$ (for which we can prove the existence of the relevant
subcomplex) to the one based on $d^{\tilde\nabla}$ (for which we have
the interpretation in terms of infinitesimal deformations).

\end{document}